\documentclass[a4paper,colorlinks,citecolor=brown
]{article}

\usepackage[utf8]{inputenc}
\usepackage[style=ieee, isbn=false, url=false,
 doi=false,sorting=nyt,backend=biber]{biblatex}
\title{Scattering diagrams in mirror symmetry}
\author{Veronica Fantini}
\date{\today}

\usepackage{tabularx}
\usepackage{amsthm}
\usepackage{amsmath}
\usepackage{amsfonts}
\newtheorem{definition}{Definition}
\newtheorem{theorem}{Theorem}
\newtheorem{remark}{Remark}

\usepackage[dvipsnames]{xcolor}
\usepackage{tikz-cd}
\usetikzlibrary{decorations.pathmorphing}

\tikzset{snake it/.style={decorate, decoration=snake}}
\usepackage{graphicx}
\usepackage{amssymb}
\usepackage[colorlinks=true,breaklinks=true,linkcolor=ForestGreen]{hyperref}

\def\C{\mathbb C}
\def\R{\mathbb R}
\def\Z{\mathbb Z}
\def\O{\mathcal O}
\DeclareMathOperator{\Der}{Der}

\DeclareMathOperator{\Spec}{Spec\,}
\DeclareMathOperator{\Hom}{Hom}
\DeclareMathOperator{\Aut}{Aut}

\begin{document}
\maketitle

\begin{abstract}
Since the pioneering work of Kontsevich and Soibelman \cite{KS-affine}, scattering diagrams have started playing an important role in mirror symmetry, in particular in the study of the \textit{reconstruction problem}. This paper aims at introducing the main ideas on the subject describing the role of scattering diagrams in relation to the SYZ conjecture and the HMS conjecture. \end{abstract}

\tableofcontents

\section*{Introduction}
\label{sec:intro}
\addcontentsline{toc}{section}{\nameref{sec:intro}}

 First conjectured by string theorists, mirror symmetry has been largely studied and it remains an active field of research in geometry and physics. 

 In its early days, a great effort was devoted to rigorous formulating mirror symmetry leading to the Strominger--Yau--Zaslow conjecture \cite{SYZ}, the Kontsevich's Homological Mirror Symmetry conjecture \cite{KS-HMS}, the Gross--Siebert program \cite{GS1,GS2,GS07}, and others even beyond Calabi--Yau varieties. 

In addition, one of the first problems addressed by mathematicians was to understand how to reconstruct the mirror of a given family of Calabi--Yau varieties. For elliptic curves, Polishchuk and Zaslow proved that elliptic curves are self-mirror dual meaning the bounded derived category of coherent sheaves of an elliptic curve is equivalent to a suitably modified Fukaya category of the mirror elliptic curve. Already for K3 surfaces, the construction of the mirror is more involved and it was first computed by Kontsevich and Soilbelman for an analytic K3 over a non--archimedean field, using combinatorial structures later called \textit{scattering diagrams} \cite{KS-affine}.
As first conjectured by Strominger, Yau, and Zaslow, mirror pairs of Calabi--Yau varieties come in families, and in a suitable limit, they admit dual torus fibrations. Then, according to Kontsevich--Soibelman, scattering diagrams consist of a collection of walls in the base of the fibration decorated with automorphisms. Furthermore, they must encode enumerative geometric data of the mirror (\textit{quantum corrections}).  

The aim of this paper is to give a gentle introduction to scattering diagrams regarded as one the essential tools to study the \textit{reconstruction problem} in mirror symmetry. What really characterizes scattering diagrams is the choice of \textit{automoprhisms} that decorate the walls of the diagram. Depending on the problem, different choices of groups are possible: for instance for log Calabi--Yau surfaces either the tropical vertex group $\mathbb{V}$ of Gross, Pandharipande, and Siebert \cite{GPS10} or the quantum tropical vertex $\hat{\mathbb{V}}$ of Bousseau \cite{bou_quantum} can be chosen. With the first choice, the quantum corrections are expressed in terms of genus-zero Gromov--Witten invariants, while with the second choice higher genus Gromov--Witten invariants contribute. \footnote{In the latter case, the mirror is the \textit{deformation quantization} of the mirror constructed by Gross--Hacking--Keel using scattering diagrams in $\mathbb{V}$ \cite{GHK}. The construction developed in \cite{bou_mirror} can be considered as a generalization of the reconstruction problem, that goes beyond the one of the Strominger--Yau--Zaslow mirror symmetry.} Recently, the author has introduced an extension of $\mathbb{V}$, the so-called extended tropical vertex group $\tilde{\mathbb{V}}$. Morally, it is expected that scattering diagrams in $\tilde{\mathbb{V}}$ should prescribe quantum corrections to reconstruct the mirror of a log Calabi--Yau surface together with a holomorphic vector bundle over it.            

The paper is divided into two main sections: in the first one (Section \ref{sec:mirror symmetry}) we present a general overview of the mirror symmetry conjectures and of the \textit{reconstruction problem}. In particular, we are going to present the approach of Kontsevich--Soibelman through non--archimedean geometry, the Gross--Siebert approach through logarithmic geometry, and the analytic approach of Fukaya.
The second part (Section \ref{sec:scattering}) is then devoted to scattering diagrams, including the examples of $\mathbb{V}, \hat{\mathbb{V}}$ and $ \Tilde{\mathbb{V}}$.

\subsection*{Acknowledgement}
The author wishes to thank professors Jacopo Stoppa and Nicol\'o Sibilla for their courses respectively on \textit{GHK construction} and on an \textit{introduction to mirror symmetry} held at SISSA in Spring 2021. 
We are also thankful to the organizers of the "First UMI meeting for PhD" held in Padova in May 2022 for the opportunity to present part of this work. Finally, we thank the anonymous referee for his precious comments and corrections which improved and enriched the exposed material.  

\section{Mirror symmetry}\label{sec:mirror symmetry}

\begin{figure}[ht]
    \centering
    \begin{tikzpicture}
    \node at (0,5.5) {$X$};
        \node at (0,4) {1};
        \node[blue] at (1,4) {$h^{2,1}$};
        \node at (2,4) {$h^{1,2}$};
        \node at (3,4) {1};
        \node at (0.5,4.5) {0};
        \node[red] at (1.5,4.5) {$h^{1,1}$};
        \node at (2.5,4.5) {0};
        \node at (0.5,3.5) {0};
        \node at (1.5,3.5) {$h^{1,1}$};
        \node at (2.5,3.5) {0};
        \node at (1,5) {0};
        \node at (2,5) {0};
        \node at (1,3) {0};
        \node at (2,3) {0};
        \node at (1.5,5.5) {1};
        \node at (1.5,2.5) {1};
        \node at (6,5.5) {$\check{X}$};
        \node at (6,4) {1};
        \node[red] at (7,4) {$h^{2,1}$};
        \node at (8,4) {$h^{1,2}$};
        \node at (9,4) {1};
        \node at (6.5,4.5) {0};
        \node[blue] at (7.5,4.5) {$h^{1,1}$};
        \node at (8.5,4.5) {0};
        \node at (6.5,3.5) {0};
        \node at (7.5,3.5) {$h^{1,1}$};
        \node at (8.5,3.5) {0};
        \node at (7,5) {0};
        \node at (8,5) {0};
        \node at (7,3) {0};
        \node at (8,3) {0};
        \node at (7.5,5.5) {1};
        \node at (7.5,2.5) {1};
    \end{tikzpicture}
    \caption{The Hodge diamonds of a mirror pair of Calabi--Yau 3-fold are symmetric with respect to the diagonal axis as $h^{1,1}(X)= h^{2,1}(\check{X})$ and $h^{2,1}(X)=h^{1,1}(\check{X})$. }
    \label{fig:diamonds}
\end{figure}
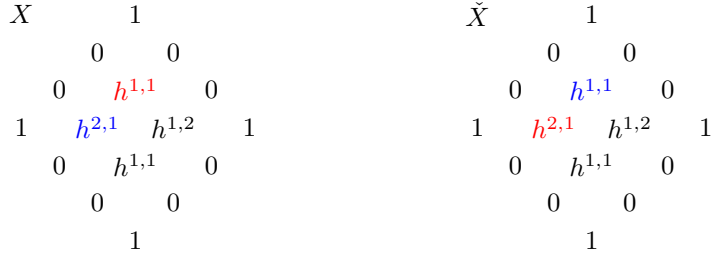

Mirror symmetry was first conjectured by string theorists\footnote{See \cite{candelas1,greene1990duality}; for more general introductive material see for example the following books \cite{morrison1992essays,CoxKatz,Dbranes,MS}. Other references will be mentioned later on.}, who observed that \textit{mirror pairs} of Calabi--Yau\footnote{Calabi--Yau manifolds are complex K\"ahler manifolds which admit Ricci flat K\"ahler metric. We denote by $\omega$ the K\"ahler form and by $J$ the complex structure. We refer to \cite{gross2012calabi} for a nice introduction on Calabi--Yau geometry with an eye toward mirror symmetry (in particular in Chapter 2).} 3--folds $(X,\omega,J)$ and $(\check{X},\check{J},\check{\omega})$ give the same theory. As a first definition of \textit{mirror pairs} we may say that $X$ and $\check{X}$ have mirror Hodge diamonds as in Figure \ref{fig:diamonds}, namely $h^{1,1}(X)=h^{2,1}(\check{X})$ and $h^{2,1}(X)=h^{1,1}(\check{X})$. 

This symmetry underlies a duality between the complex moduli and the complexified K\"ahler moduli between the mirror pairs: indeed on the one hand the deformations of the complex structure of $X$ are parametrized by $H^1(X,T_X)$ and by Serre duality and Dolbeault's theorem we have the following chain of equivalences
\begin{equation*}
    H^1(X,T_X)\cong H^{1}(X,\Omega_X^2)^{\vee}\cong H^{2,1}(X)^{\vee}
\end{equation*}
On the other hand, the K\"ahler moduli is the cone of K\"ahler classes $\omega\in H^{1,1}(X,\C)\cap H^2(X,\R)$ associated to a fixed K\"ahler metric on $X$
\begin{equation*}
    \mathcal{K}_X:=\lbrace [\omega]\in H^{1,1} \vert [\omega] \text{ is represented by a K\"ahler form}\rbrace\underset{open}{\subset} H^{1,1}(X,\R)
\end{equation*}
Since the complex moduli is a complex manifold, it is necessary to consider the complexified K\"ahler moduli defined as \[\mathcal{K}_{X}^{\C}:=(H^2(X,\R)\oplus i\mathcal{K}_X)/H^2(X,\Z).\] 
Recall that if $X$ is a smooth projective variety, the number $h^{p,q}(X)$ are the dimension of complex vector space $H^{p,q}(X)$. Therefore, having mirror Hodge diamonds reflects the following duality between moduli 
\begin{align*}
    \dim H^{1}(X,T_X)=h^{2,1}(X)& = h^{1,1}(\check{X})=\dim  H^2(\check{X},\C)=\dim {K}_{\check{X}}^{\C} \\
    \dim H^{1}(\check{X},T_{\check{X}})=h^{2,1}(\check{X})& = h^{1,1}(X)=\dim H^2(X,\C)=\dim \mathcal{K}_{X}^{\C} 
\end{align*}

The breakthrough, attracting the interest of the geometry community, was the paper by Candelas, De la Ossa, Green, and Parkers \cite{COGP} who were able to compute the number of rational curves (up to degree 9) on a quintic 3-fold by studying solutions of Picard--Fuchs equations. The number of rational curves (also known as Gromov--Witten invariants) are invariants that only depend on the symplectic form $\omega$, and with the algebraic-geometric techniques available at that time only the number of curves up to degree 3 was known. 
While the Picard--Fuchs equations are equations for periods of $\check{X}$ with respect to the holomorphic volume form of $(\check{X},\check{J})$. Hence we see again mirror symmetry as a duality between complex $(\check{X},\check{J})$ and the symplectic $(X,\omega)$ (and vice-versa between $(\check{X},\check{\omega})$ and $(X,J)$).

So far we have seen one of the first aspects of mirror symmetry for Calabi--Yau 3-folds, but later on, different people worked on a geometric formulation of the mirror symmetry conjecture. Among them we will recover the Strominger--Yau--Zaslow (SYZ) conjecture \cite{SYZ} and the Homological Mirror Symmetry of Kontsevich \cite{K95-HMS}. 

In addition, geometers started looking for examples of mirror pairs (beyond the 3-folds) by studying the so-called \textit{reconstruction problem}, namely how to build the mirror of a given Calabi--Yau. We will recover some of the main results addressing the reconstruction problem: Kontsevich--Soibelman mirror symmetry via non--archimedean geometry, the Gross--Siebert program via logarithmic geometry, and Fukaya's mirror symmetry via multi-valued Morse theory. As discussed in the following sections, their common aspect can be traced back to the use of \textit{scattering diagrams}. 

More recently, other progresses have been pursued in the study of the reconstruction problem from different perspectives and using techniques that go beyond scattering diagrams. We will recall some of them in the exposition but without giving many details.

\subsection{Strominger--Yau--Zaslow conjecture}

The Strominger--Yau--Zaslow conjecture (SYZ) claims that mirror symmetry is T-duality, namely that pairs of Calabi--Yau 3-folds $X$ and $\check{X}$ carry special Lagrangian torus fibrations (with possible singular fibers) $\pi\colon X\to B$ and $\check{p}\colon\check{X}\to B$. Away from singularities, T-duality acts on the smooth fibers as a map from a Lagrangian fiber $X_b=\pi^{-1}(b)$ to the moduli of unitary flat connections on a special Lagrangian fiber of the mirror $H^1(\check{X}_b,\mathbb{R}/\mathbb{Z})$. However, due to the presence of singularities, it is necessary to add corrections \footnote{In the 1996 paper \cite{SYZ} \textit{quantum corrections} are needed to construct the Ricci-flat K\"ahler metric on the mirror. Although the metric approach seemed originally out of reach in higher dimensions (handling quantum corrections seemed also difficult due to the fact that already in 3 dimensions the critical locus of Lagrangians fibrations is not of codimension 2 \cite{joyce2000singularities}), new progress has been made by Collins, Jacobs and Lin \cite{CJL21,CJL20} but mostly in 2 dimensions. }. These should encode enumerative geometric data such as counting holomorphic discs. 

We illustrate this duality in the following toy model example, where we avoid singularities. 

\subsubsection*{Toy model I} Let $B_0$ be an integral affine smooth manifold (i.e. with transition functions in $GL(n,\Z)\ltimes \R^n$). Let $\Lambda$ be an integer lattice subbundle of the tangent bundle of $B_0$,  $\Lambda\subset TB_0$, and its dual lattice $\Lambda^*=\Hom(\Lambda,\Z)\subset T^*B_0$. Then we define two torus fibrations 

\begin{align*}
    X:=T^*B_0/\Lambda^*\to B_0 &\qquad  \check{X}:=TB_0/\Lambda\to B_0
\end{align*}

such that $X$ comes with a natural symplectic form $\omega$ and $\check{X}$ with a natural complex structure $J$; in local affine coordinates $(x^i,y^i)$ 

\begin{align*}
    \omega=\sum_{i} dx^i\wedge dy^i & \qquad \check{J}\colon \frac{\partial}{\partial x_j}\to -\frac{\partial}{\partial y_j}. 
\end{align*}

Then $(X,\omega)$ and $(\check{X},\check{J})$ are mirror pairs that admit dual torus fibrations and such that the symplectic coordinates of $X$ are the complex coordinates of $\check{X}$. This is only the half of mirror symmetry; if $B_0$ admits a Riemannian metric $g$ of hessian type, i.e. there exists a function $h$ such that in local coordinates \[g_{ij}(x)=\frac{\partial^2 h}{\partial x_i\partial x_j} \] then $\check{x}_i:=\frac{\partial h}{\partial x_i}$ is a complex coordinate on $X$ with respect to the complex structure $J\left( \frac{\partial}{\partial x_i}\right)=\sum_{j} g_{ij}(x)\frac{\partial}{\partial y_j}$. Indeed,

\begin{align*}
    J\left(\frac{\partial}{\partial y_i}\right)=-\sum_{j}g^{ij}(x)\frac{\partial}{\partial x_j}= -\sum_{j}\frac{\partial^2 h}{\partial x^i\partial x^j}\frac{\partial}{\partial x_j}=-\frac{\partial}{\partial \check{x}_i}.
\end{align*}
The change of coordinates $\check{x}_i:=\frac{\partial h}{\partial x_i}$ is a Legendre transform. Furthermore, assuming $\det \left(\frac{\partial^2 h}{\partial x_i\partial x_j}\right)=1 $, $g$ gives a K\"ahler metric on $X$ compatible with $J$ and $\omega$. Conversely, let $\check{g}=\sum_j d\check{x}^j d\check{x}^j$ be another metric on $B_0$: we define \[\check{\omega}\left(\frac{\partial}{\partial \check{x}_i}, \check{J}\left(\frac{\partial}{\partial \check{x}_j}\right)\right):=\check{g}\left(\frac{\partial}{\partial \check{x}_i},\frac{\partial}{\partial \check{x}_j}\right) \]
 to be the compatible K\"ahler form with respect to $\check{g},\check{J}$. Then $\check{\omega}=\sum_j g_{ij}(x)d\check{x}^i\wedge dy^j$ hence $\check{x}^i,y^i$ are symplectic coordinates on $\check{X}$. Schematically, we represent mirror symmetry for $X$ and $\check{X}$ as follows 
 
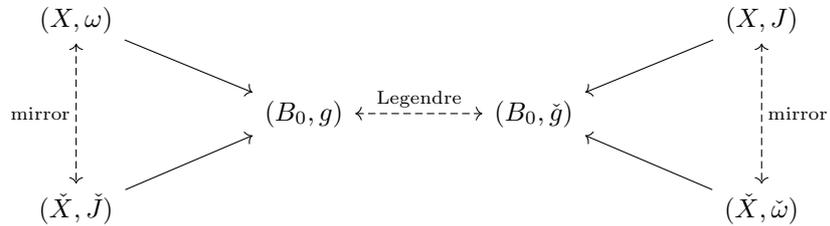
\vspace{3mm}

 \begin{tikzcd}
     (X,\omega)\arrow[dd,dashed,<->,swap,"\text{mirror}"]\arrow[drr] & &  & & & &(X,J)\arrow[dd,dashed,<->,"\text{mirror}"]\arrow[dll]  \\
     & & (B_0,g)\arrow[rr,dashed,<->,"\text{Legendre}"]& & (B_0,\check{g})& &\\
   (\check{X},\check{J})\arrow[urr]&  &  & & &   &(\check{X},\check{\omega})\arrow[ull] 
 \end{tikzcd}

\subsection{Kontsevich--Soibelman approach via non--archimedean geometry}

Following the proposal of the SYZ conjecture, in \cite{KS-HMS} Kontsevich and Soibelman propose to replace the dual SYZ torus fibration with a non-archimedean torus fibration. In this framework, they address the reconstruction problem of K3 surfaces, whose mirror will be constructed as an analytic space over a non--archimedean field \cite{KS-affine}.\footnote{This approach successfully applies in \cite{tu14} where the author constructs the mirror of a singular Lagrangian torus fibration as a rigid analytic space obtained by gluing local affine pieces with transition functions encoding counting pseudo-holomorphic discs. However, the mirror space is constructed on the smooth locus of the base.} Their construction is built on some additional structures on the base of the SYZ fibration, namely the so-called \textit{scattering diagrams}. \footnote{In fact the name \textit{scattering diagrams} does not appear in the original work of Kontsevich and Soibelman, but it was later introduced by Gross and Siebert in \cite{GS07}. In \cite{KS-affine} Section $9$, the authors rather talk about \textit{lines} on the base of the fibration as gradient lines satisfying certain axioms. In the following, we will refer to \textit{rays} as the image of gradient lines under a local homeomorphism from $U\subset B$ to $\R^2$ (see \cite[Axiom 6]{KS-affine}).}

As suggested by SYZ conjecture, mirror pairs of Calabi--Yau admit locally dual torus fibrations; these arise naturally in both rigid analytic complex and symplectic geometry. On the one hand, they consist of a smooth map $\pi\colon X\to B$ from  a variety $X$ defined over non--archimedean field to a CW complex $B$. On the other hand, they are integrable systems. In addition, in both cases, the base manifold $B$ carries an integral affine structure with singularities. Conjecturally, integral affine manifolds with singularities arise as the Gromov--Hausdorff limits of families of Calabi--Yau (see \cite[Conjecture 1]{KS-affine} and \cite[Conjecture 6.2]{gross-wilson}, in fact for K3 surfaces the Conjecture was proved by Gross and Wilson \cite[Theorem 6.4]{gross-wilson}). Hence, away from singularities, it should be possible to locally reconstruct the torus fibration from the base manifold, more precisely from its integral affine structure. However, due to the presence of singularities, it is necessary to consider additional structures in order to construct globally the fibration. Here is where the scattering diagrams come into play.  


We now discuss the example of K3 surfaces, following \cite{KS-affine}. From the symplectic side, elliptically fibered K3 surface $(X,\Omega_X)$ gives a complex integrable system $(X,\C\mathbb{P}^1,\omega)$: a holomorphic fibration $\pi\colon X\to \C\mathbb{P}^1$ whose fibers are elliptic curves, $\omega=\mathrm{Re } \,\Omega_X$ with $\Omega_X$ the holomorphic $2$-form on $X$ and $B=\C\mathbb{P}^1$ has 24 singular fibers. The integral affine structure on $B$ is induced from the one of $\R^2$ with nontrivial monodromy around the singularities. Then Kontsevich and Soibelman construct the mirror $\check{X}$ as an analytic space over $K=\C(\!(t)\!)$. 
We recall the heuristic of the main steps of their construction and we refer to \cite{KS-affine} for details:
\begin{itemize}
   \item on the smooth locus $B_0\subset B$ the mirror is locally isomorphic to the dual torus fibration: the group $SL(2,\Z)\ltimes (K^\times)^2$ (induced by the integral affine structure of $B_0$) acts by automorphisms of the fibration $\pi_{can}\colon(\mathbb{G}_m^{an})^2\to \R^2$. 
    \item in a neighborhood of the singularity $U_s\subset B$, the local sheaf of functions $\O^{sing}$ is defined by a special choice of gluing automorphisms.
    \item in order to glue the local model of singularities to the canonical one, they introduce a new structure: consider a ray emanating from the singularity together with an automorphism and deform the sheaf of functions across the ray (as in Figure \ref{fig:my_label1}). 
    \begin{figure}[ht]
        \centering
        \begin{tikzpicture}
            \node[red, thick] at (0,0) {$\times$};
            \draw (0,0)--(3,1.5);
            \draw[gray] (0,0) circle (0.5);
            \draw[gray] (0.8,0.4) circle (0.5);
            \node[font=\small,below,gray] at (1.5,0.5) {$U$}; 
             \node[font=\small,below,gray] at (0,-0.5) {$U_s$}; 
        \end{tikzpicture}
        \caption{Ray emanated from the singularity. The action of the automorphism that decorates the ray modifies the canonical sheaf of function over $U$ in the sheaf of functions $\O^{mod}$ }
        \label{fig:my_label1}
    \end{figure}
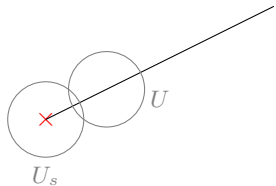
    \item assuming that two rays meet transversely, the action of the two automorphisms $g_1,g_2$ must be consistent, i.e. the modified sheaves $\O^{mod}_1, \O^{mod}_2$ must glue together in a neighborhood of the intersection point. If $g_1,g_2$ commute, then $\O^{mod}_1, \O^{mod}_2$ glue well. On the contrary, if $g_1,g_2$ do not commute, then the deformation of $\O^{can}$ must receive other (possibly infinitely many) corrections by $g_1,g_2,...,g_k,...$ with the property that \[g_1\cdot g_2\cdot ... \cdot g_k \cdot ...\cdot g_2^{-1}\cdot g_1^{-1}=\mathrm{id}_{\mathbb{G}}.\] The latter property assures that the monodromy of the fiber above the intersection point remains trivial. Locally, the set of lines and rays meeting at a smooth point, together with the automorphisms $g_k$ is called a \textit{scattering diagram} (see Figure \ref{fig:my_label2}). 
    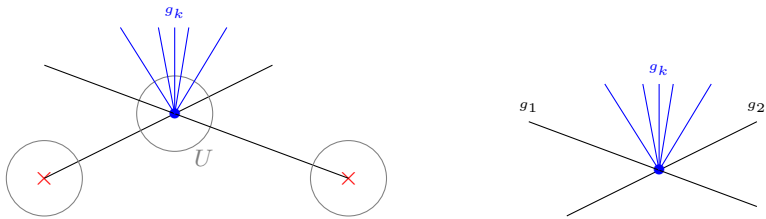
\begin{figure}[ht]
        \centering
        \begin{tikzpicture}
            \node[red, thick] at (0,0) {$\times$};
            \draw (0,0)--(3,1.5);
            \draw[gray] (0,0) circle (0.5);
            \draw[gray] (12/7,6/7) circle (0.5);
            \node[thick, blue] at (12/7,6/7) {$\bullet$};
             \node[red, thick] at (4,0) {$\times$};
             \draw (4,0)--(0,1.5);
            \draw[gray] (4,0) circle (0.5);
            \draw[blue] (12/7,6/7) --(12/7,2);
            \draw[blue] (12/7,6/7) --(1.5,2);
            \draw[blue] (12/7,6/7) --(2.4,2);
            \draw[blue] (12/7,6/7) --(1.9,2);
            \draw[blue] (12/7,6/7) --(1,2);
            \node[blue,above, font=\tiny] at (12/7,2) {$g_k$};
            \node[font=\small,below,gray] at (2.1,0.5) {$U$};
                   \end{tikzpicture}
                   \hspace{1.5cm}%
                   \begin{tikzpicture}
            \draw (0.5,0.25)--(3,1.5);
            \node[thick, blue] at (12/7,6/7) {$\bullet$};
            \draw (3,3/8)--(0,1.5);
            \draw[blue] (12/7,6/7) --(12/7,2);
            \draw[blue] (12/7,6/7) --(1.5,2);
            \draw[blue] (12/7,6/7) --(2.4,2);
            \draw[blue] (12/7,6/7) --(1.9,2);
            \draw[blue] (12/7,6/7) --(1,2);
            \node[blue,above, font=\tiny] at (12/7,2) {$g_k$};
              \node[above, font=\tiny] at (3,1.5) {$g_2$};
                \node[above, font=\tiny] at (0,1.5) {$g_1$};
                   \end{tikzpicture}
        \caption{On the left is the local picture of the rays meeting at a generic point and the new rays that are needed to deform the canonical sheaf over $U$. In fact, we draw actual rays in $\R^2$ (and not wiggly lines in $B$ as in Figure 6 of \cite{KS-affine}) under the assumption that there is a local homeomorphism from $U\subset B$ to $\R^2$ (see \cite[Axiom 6]{KS-affine}). On the right, we zoom in the neighborhood of the intersection point drawing a local scattering diagram, where the location of the singularity is neglected.}
        \label{fig:my_label2}
    \end{figure}
    \item finally they construct a nowhere vanishing top degree analytic form $\Omega$, which extends over the modified sheaves without modifying the $K$-affine structure of $B$. 
\end{itemize}

We conclude this section with a simple example: the elliptic curve. 

\subsubsection*{Example: elliptic curve I}\label{elliptic curve}

Let $\pi\colon X\to B$ be an affine map of tori with $X:=\R^2/\Lambda$ and $\Lambda\cong\Z^2$, $B=\R/\Lambda'$ and $\Lambda'\cong\Z$ (in particular we assume $\Lambda'=\alpha\Z$ for $\alpha\in\R\setminus\{0\}$). The base $B$ is smooth and it has a natural affine structure induced by $\Lambda'$.Then, the dual torus fibration $\check{\pi}\colon \check{X}\to B$ can be easily built from the integral affine structure of $B$.

Let $K=\C(\!(q)\!)$ be the ring of Laurent polynomial in $q$, and let $\rho\colon\Z\to K^\times$ be a group homomorphism such that the image of the composition $\mathit{val}\circ\rho\colon\Z\to\R$ is the lattice $\Lambda'\subset\R$. Then $q^{\Lambda '}=\{q^{n}\}_{n\in\Lambda'}$ acts on $\C^\times$ by translations and the quotient is a K-analytic space called the Tate curve $\check{X}:=\C^{\times}/q^{\Lambda'}$. The dual fibration is $\check{\pi}\colon \check{X}\to B$. 

Since every elliptic curve over $\C$ is isomorphic to $ \C/\Lambda_\tau$ with $\Lambda_\tau=\Z+\tau\Z$ and $\tau$ in the upper half-plane, the map $\mathbf{e}\colon z\to e^{2\pi iz}$ induces an isomorphism $\C/\Lambda_{\tau}\cong\check{X}$ with $q=\mathbf{e}(\frac{\tau}{\alpha})$.  



\subsection{Gross--Siebert approach via logarithmic geometry}

Scattering diagrams appear also in the Gross--Siebert program \cite{GS1,GS2,GS07} \footnote{The reader may find useful the following readings \cite{GSintrinsic,GS-review}.} which solves the reconstruction problem through logarithmic geometry and in all dimensions.\footnote{Although running in a very different direction, the recent work of Groman and Varulgones \cite{groman2021locality,groman2022closed} addresses the reconstruction problem in the spirit of the Gross--Sibert program.} In particular, they replace the set of \textit{rays} of \cite{KS-affine} (which is limited to dimension $2$) by a codimension one polyhedral complex of \textit{walls}. 

In particular, an interesting class of examples is given by log Calabi--Yau pairs. In two dimensions the reconstruction of the mirror of log Calabi--Yau surfaces has been studied in \cite{GHK,bou_mirror} (among which some toric examples are computed in \cite{barrott}) and of non-compact 3-folds as local $\mathbb{P}^2$ (see \cite[Examples 5.1-5.2]{GSinvitation} and \cite{GSlocal}). In higher dimension, some examples of computing the mirror of log Calabi--Yau pairs are given in \cite{arguz-higher-logCY}. We review below the construction of \cite{GHK}, which reinterpret the construction of the Gross--Siebert program in terms of curve counting (crucially using \cite{GPS10}), \textit{broken lines} and generalized theta functions (introduced in \cite{Gross-P2}). 

Let $(Y,D=D_1+...+D_N)$ where $Y$ is a smooth projective surface and $D\in |-K_Y|$ is a cycle of rational curves, then $X=Y\setminus D$ is a log Calabi--Yau surface and it admits a holomorphic symplectic form $\Omega_X$ with simple poles along the divisor $D$. So $(X,\Omega_X)$ is non compact Calabi--Yau. The toy model examples of log CY are toric surfaces and their toric boundary such as $\mathbb{P}^2, \mathbb{P}^1\times\mathbb{P}^1$ and $dP_7$ (which is the blow-up of $\mathbb{P}^2$ in two distinct points). More generally, for a given log CY surface $(Y,D)$ there exists a toric blow-up $(\tilde{Y},\tilde{D})$ which has a toric model $(\bar{Y},\bar{D})$ (see the diagram below). Examples of non-toric log CY are the del Pezzo surfaces $dP_5, dP_4,dP_3, dP_2$. 

\begin{figure}[ht]
    \centering
    \begin{tikzcd}
        &  (\tilde{Y},\tilde{D})\arrow[rd,"p"]\arrow[dl, dashed, swap, "\pi"] & \\
        (Y,D) & & (\bar{Y},\bar{D})
    \end{tikzcd}
    \caption{The map $\pi\colon \tilde{Y}\to Y$ is the toric blow-up, namely it is a birational map such that $(\pi^{-1}(D))^{\text{red}}$ is an anticanonical cycle of rational curves. The map $p\colon\tilde{Y}\to \bar{Y}$ is the blow-up of $\bar{Y}$ at some distinct points $x_{ij}\in \bar{D}_i$ and $\tilde{D}$ is the proper transform of $\bar{D}$. }
    \label{fig:my_label}
\end{figure}
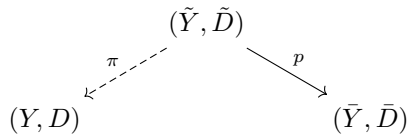

Schematically, we summarize the Gross--Siebert program as follows: it takes place at two different levels, namely the mirror map between families of Calabi--Yau $\mathcal{X}\longleftrightarrow \check{\mathcal{X}}$ is already incarnated at the tropical level through a discrete Legendre transform $B\longleftrightarrow \check{B}$.    

\begin{center}
\begin{tikzcd}
\mathcal{X}\arrow[rr, dashed,<->,"\textit{mirror}"]\arrow[d,dashed] & & \check{\mathcal{X}} \\
    (B,\varphi)\arrow[rr,<->,"\textit{Legendre}"] & & (\check{B},\check{\varphi})\arrow[u,dashed]
\end{tikzcd}
\end{center}

Most important, in order to make sense of the vertical arrows of the diagram it is necessary to use logarithmic geometry. Furthermore, log geometry plays a fundamental role in proving the \textit{numerical} correspondence between Gromov--Witten invariants and Hodge periods of the mirror. 

On the one hand, Gross and Siebert define a discrete Legendre transform as \textit{mirror} map between $(B,\mathcal{P},\varphi)$ and $(\check{B},\check{\mathcal{P}},\check{\varphi})$ where $B$ is a tropical affine manifold with singularities, $\mathcal{P}$ is a polyhedral decomposition of $B$ and $\varphi$ is a convex multi-valued piecewise linear function on $B$. On the other hand, they show how to build a smooth family $\check{\mathcal{X}}$ from the data of $(\check{B},\check{\mathcal{P}})$ and of a gluing map $\check{s}$ and vice-versa how to recover the data of $(B,\mathcal{P},\varphi)$ and $s$ from the family $\mathcal{X}$. In particular, from $(\check{B},\check{\mathcal{P}},\check{s})$ they construct a log scheme $\check{X}_0(\check{B},\check{\mathcal{P}},\check{s})^{\dagger}$ which can be smoothed into a formal family $\check{\mathfrak{X}}$. Furthermore, if the family admits a polarization, then there will be an underlying scheme $\check{\mathcal{X}}$. Conversely, starting with a family $\mathcal{X}$ of Calabi--Yau varieties which admit a \textit{toric degeneration} i.e. whose central fiber $X_0^{\dagger}$ is a union of toric varieties and the map to the base is log smooth away from a singular locus, they recover $({B},{\mathcal{P}})$ and a gluing data $s$ such that $X_0=X_0(B,\mathcal{P},s)$. If $\mathcal{X}$ is polarized then $B$ comes equipped with a convex multi-valued piecewise linear function $\varphi$. Pictorially, we may sketch the Gross--Siebert construction as follows 

\begin{center}
\begin{tikzcd}
(\mathcal{X}, \mathcal{L}_\blacktriangle)\arrow[rr, dashed,<->,"\textit{mirror}"]\arrow[d,swap, dashed, "\text{toric degen.}"] & & (\check{\mathcal{X}},\check{\mathcal{L}}) \\
(X_0^{\dagger},s)\arrow[d,dashed]& & (\check{X}_0^\dagger,\check{s}_{\bullet})\arrow[u,dashed]\\
    (B,\mathcal{P},\varphi_{\blacktriangle})\arrow[rr,<->,"\textit{Legendre}"] & & (\check{B},\check{\mathcal{P}}, \check{\varphi})\arrow[u,dashed,swap,"\bullet"]
\end{tikzcd}
\end{center}

where $\blacktriangle$ denotes a choice of polarization and $\bullet$ a choice of gluing $\check{s}$. 

Although, in full generality the Gross--Siebert program is quite technical, in the case when $X=Y\setminus D$ is a toric log CY surface and $D=D_1+...+D_N$, it simplifies as sketched below: let $\Sigma$ be the fan of the toric surface $Y$ and $\Delta_{(Y,D)}$ be the dual intersection complex 
\begin{itemize}
    \item the base $B$ is the topological space underlying $\Delta_{(Y,D)}$,  and it admits an integral affine structure that depends on the intersection matrix of $D$.   
    \item the base of mirror family is $\Spec \C[\mathrm{NE}(Y)]$ where $\mathrm{NE}(Y)$ is the cone of effective curve classes (which is finitely generated by toric divisors)
    \item the construction of the mirror family $\check{\mathcal{X}}\to\Spec\C[\mathrm{NE}(Y)]$ mimics toric Mumford's degenerations \cite{GSinvitation}, and it is obtained as a smoothing of the $N$-vertex $\mathbb{V}_N$\footnote{The $N$-vertex $\mathbb{V}_N:=\mathbb{A}^2_{x_1,x_2}\cup ...\cup\mathbb{A}^2_{x_{N-1},x_N}\cup \mathbb{A}^2_{x_N,x_1}$ is defined as the union of $N$ coordinates hyperplanes.}. In particular, the smooth fibers are algebraic tori $\C^*\times\C^*$ (and they are indeed \textit{mirror} of $X$). 
\end{itemize}

If $X=Y\setminus D$ is a log CY surface (non-toric) and $D=D_1+...+D_N$, the previous construction has to be modified. First of all, recall that $\mathrm{NE}_{\geq 0}(Y)$ is not always a finitely generated monoid (for instance if the adjacency matrix is not negetive semi-definite, then $\mathrm{NE}_{\geq 0}(Y)$ is finitely generated). So more generally, one has to consider $\sigma_P\subset A_1(Y,\R)$ a strictly convex rational polyhedral cone containing $\mathrm{NE}(Y)_{\R}$ and its associated monoid $P=\sigma_P\cap  A_1(Y,\Z)$. Then, by modifying Mumford degeneration, the mirror family is constructed as a smoothing of $\mathbb{V}_N^0:=\mathbb{V}_N\setminus\{0\}$. However, the smooth fibers must receive corrections from consistent scattering diagrams in order to glue them together with the singular fiber at the origin (details in \cite{GHK}). The outcome of this procedure gives a family $\check{X}_I\to \Spec R/I$ where $R:=\C[P]$ and $I\subset R$ is a monomial ideal. Most importantly,   $\check{X}_I$ is defined explicitly as the spectrum of an algebra on a vector space with a canonical basis (the basis of \textit{theta functions}) and a product defined in terms of the Gromov--Witten theory of $(Y, D)$. Finally, the resulting family $\check{\mathcal{X}}\to\mathrm{Spf} \hat{R}$ is a formal family and $\hat{R}$ is the completion of $R$ with respect to the ideal $\sqrt{I}$.

\begin{remark}
    After the development of punctured Gromov--Witten invariants \cite{abramovich2020punctured}, in \cite{GS19} the authors address the problem of reconstructing the mirror of a log Calabi--Yau pair by constructing directly the ring of theta functions; this program goes under name of {\em intrinsic mirror symmetry}\footnote{An equivalent construction using non--archimedean methods (on the symplectic side) is due to Keel and Yu \cite{KY19} when $X$ contains an algebraic torus.}. It is remarkable, as proved in \cite{GS22}, that intrinsic mirror symmetry for log Calabi--Yau pair fits into the Gross--Siebert program, namely it constructs the same mirror as the one obtained using scattering diagrams in the spirit of the GHK construction we briefly discussed (but in all dimensions). 
\end{remark}

\begin{remark}
    In the GHK construction, the role of consistent scattering diagrams is essential to obtain the formal mirror family $\check{\mathcal{X}}\to \mathrm{Spf}\hat{R}$ by smoothing the $N$-vertex $\mathbb{V}_N$. More generally, degeneration of Calabi--Yau manifolds (beyond log Calabi--Yau surfaces), gives a singular Calabi--Yau variety equipped with a natural log structure. In \cite{FFR}, using a deformation theory approach with differential graded Batalin--Vilkovisky algebras, it has been recently proved that the smoothing is possible for \textit{toroidal crossing spaces} (both $d$-semistable log smooth Calabi--Yau varieties and maximally degenerate ones are examples of \textit{toroidal crossing spaces} and their smoothing was also studied in \cite{chan21--smoothing}). 
\end{remark}

\subsection{Fukaya approach via multi-valued Morse theory:\\
closed string conjecture}

In \cite{fuk05}, the author proposed a new approach to SYZ mirror symmetry based on the asymptotics of the solutions of Maurer--Cartan equation. As already discussed, one of the main difficulties in reconstructing the mirror manifold is due to the presence of singularities on $B$, and Fukaya proposes to study the corrections via \textit{multi-valued Morse theory}. Also in this case, the new structure consists of some combinatorial data on the integral affine manifold $B$ (or of its Legendre dual $\check{B}$) which plays the role of scattering diagrams. Roughly speaking, in Fukaya's picture the rays of scattering diagrams are replaced by trivalent graphs embedded in $B$ with roots at the singular points, and the automorphisms which decorate the rays (of scattering diagrams) are replaced by gradient flow vector fields associated with multi-valued functions $f$ on $B$. 
Furthermore, Fukaya conjectures that this combinatorial structure on the one hand encodes pseudo-holomorphic discs bounding the fibers of $X$ and on the other hand, it governs infinitesimal deformations of the mirror $\check{X}$. This goes under the name of \textit{closed string conjecture} (see \cite[Conjecture 2.2]{fuk05}).  

\begin{figure}[ht]
\center
\begin{tikzpicture}
\draw[blue, thick] (0,0) rectangle (3,1.5);
\draw[magenta, thick] (4,0) rectangle (7,1.5);
\draw[magenta, thick] (4,-2) rectangle (7,-0.5);
\draw[red, thick] (8,0) rectangle (11,1.5);
\node[font=\small] at (1.5,0.9) {Pseudo-holomorphic };
\node[font=\small] at (1.5,0.6) { discs on $X$ };
\node[font=\small] at (5.5,1.1) {Multi-valued  };
\node[font=\small] at (5.5,0.8) { Morse theory };
\node[font=\small] at (5.5,0.5) {on $B$};
\node[font=\small] at (5.5,-1) {Tropical geometry };
\node[font=\small] at (5.5,-1.3) { on Legendre dual $\check{B}$ };
\node[font=\small] at (9.5,0.9) {Deformation theory  };
\node[font=\small] at (9.5,0.6) { of $\check{X}$ };
\draw[blue, <->]  (3,0.75) -- (4,0.75);
\draw[red, <->]  (7,0.75) -- (8,0.75);
\draw[magenta, <->]  (5.5,0) -- (5.5,-0.5);
\end{tikzpicture}
\end{figure}

Recently, K. Chan, N. C. Leung and Z. N. Ma proposed a new version of Fukaya's closed string conjecture by replacing multi-valued Morse theory by scattering diagrams in $B$. In \cite{CLM19}, they first study the relationship between scattering diagrams and infinitesimal deformations of $\check{X}:=TB_0/\Lambda$ locally on an open affine subset of the smooth locus $B_0\subset B$. In \cite{chan2020tropical}, the authors proved that solutions of Maurer--Cartan equation which governs deformations of $\check{X}$ encodes tropical discs counting. The global picture is then obtained in \cite{chan2022smoothing} by relating scattering diagrams and solutions of Maurer--Cartan equation which governs deformations of log Calabi--Yau variety (in the spirit of the Gross--Siebert program). Therefore all together these results give a new formulation of Fukaya's closed-string conjecture.          

\subsection{Homological Mirror Symmetry of Kontsevich: mirror symmetry as an equivalence of categories}\label{HMS}

In \cite{K95-HMS}, Kontsevich proposed a categorical formulation of mirror symmetry, the Homological Mirror Symmetry (HMS), namely an equivalence between the bounded derived category of coherent sheaves $\mathcal{D}^{b}(\check{X})$ and the Fukaya category of the mirror $\mathcal{F}(X)$. A precise statement of HMS conjecture needs actually some explanations as on the one hand $\mathcal{D}^{b}$ is a triangulated category, and on the other hand, $\mathcal{F}$ is an $A_\infty$ category. 
We are not going to give a precise statement in the general case, and we refer to \cite{Dbranes,KS-HMS,auroux2009special,seidel2008fukaya}, but we are going to present the example of the elliptic curve where the statement simplifies. 
Looking simply at the objects of the aforementioned categories, HMS predicts a duality between coherent sheaves on $\check{X}$ and Lagrangians on the mirror $X$ \footnote{ In string theory, one will refer to this duality as a duality between A-branes (Lagrangians with flat connections) and B-branes (coherent sheaves) \cite{Dbranes}.} hence, at a first glance we will enhance the SYZ toy model by considering Lagrangians as mirrors to holomorphic vector bundles.  

\subsubsection*{Toy model II}
Let $B_0$ be an integral affine smooth $1$-dimensional manifold. Let $\Lambda\subset TB_0$ be a lattice subbundle of the tangent bundle of $B_0$ and $\Lambda^*:=\Hom(\Lambda,\Z)$ be the dual lattice. 
On the one hand, we can define a Lagrangian torus fibration $\pi\colon X:=T^*B_0/\Lambda^*\to B_0$ and $r$ \textit{horizontal} Lagrangians $L_1,...,L_r$ in $X$ with unitary local systems. 
On the other hand, the dual torus fibration is $\check{\pi}\colon\check{X}:=TB_0/\Lambda\to B_0$ and the mirror of $L_1\cup L_2\cup...\cup L_r$ is a holomorphic trivial rank $r$ vector bundle $\check{E}=\O_{\check{X}}\oplus...\oplus\O_{\check{X}}\to\check{X}$. 

\begin{figure}[ht]
    \centering
    \begin{tikzpicture}
        \draw (0,0) rectangle (3,3);
        \draw[green] (0,0)--(3,0);
        \draw[red] (0,0.8)--(3,0.8);
        \draw[blue] (0,2)--(3,2);
        \node[blue,above,font=\small] at (0.7,2) {$L_r$};
        \node[green,above,font=\small] at (2.5,0) {$L_1$};
        \node[red,above,font=\small] at (1.5,0.8) {$L_2$};
        \node[font=\small] at (1.5,0) {$\sim$};
        \node[font=\small] at (1.5,3) {$\sim$};
        \node[font=\small] at (0,1.5) {$/$};
        \node[font=\small] at (3,1.5) {$/$};
        \node[font=\small] at (3.6,3) {$X\subset\R^2$};
        \node at (0,0) {$\bullet$};
        \node[below, font=\small] at (0,0) {$\text{e}$};
        \node[font=\small] at (7,0.5) {$L_1 \leadsto \O$ };
       \node[font=\small] at (7,1.5) {$(L_2, 2\pi i\beta_2 dx) \leadsto (\O, \theta_2)$ };
       \node[font=\small] at (7,2.5) {$(L_3, 2\pi i\beta_3 dx) \leadsto (\O, \theta_3)$};
    \end{tikzpicture}
    \caption{On the left-hand side the fundamental domain of the torus $X$ and the horizontal Lagrangians $L_j=\{\beta_j\}\times [0,1]$. A part from $L_1$ which has a trivial flat connection ($\beta_1=0$), the other Lagrangians come with the flat connection $2\pi i\beta_j dx$. Their mirrors as objects of $\mathcal{D}^b(\check{X})$ are the sheaf $\O_{\check{X}}$ together with a global section $\theta_j$, which characterizes (different) line bundles in $\mathrm{Pic}^0$.}
    \label{fig:horizontal_Lag}
\end{figure}
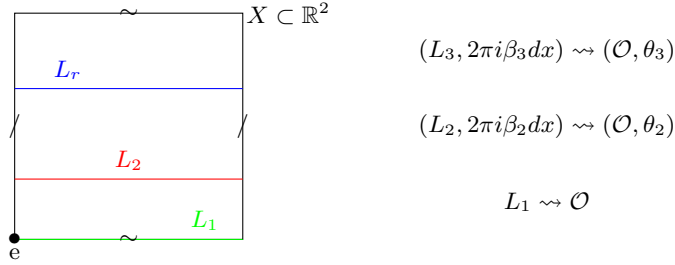

If $X=\R^2/\Z\oplus\Z$ is a symplectic torus with $\omega=dx\wedge dy $, then we can represent the horizontal Lagrangians $L_j$ as in Figure \ref{fig:horizontal_Lag}. Each Lagrangian $L_j$ comes equipped with a flat line bundle with connection $2\pi i\beta_j dx$ where $x\sim x+1$ is a local coordinate on $L_j$. 

Following \cite{PZ00}, the mirror of $X$ is the elliptic curve $\check{X}=\C/\Z+ i\Z$ and the trivial rank $r$ bundle $\check{E}\to\check{X}$ whose base of global sections is generated by theta functions $\theta_j$: \[\theta_j:=\theta[0,\beta_j](\tau,z)=\sum_{m\in\Z}\exp\{2\pi i[\tau \frac{m^2}{2}+m(z+\beta_j)]\}.\]

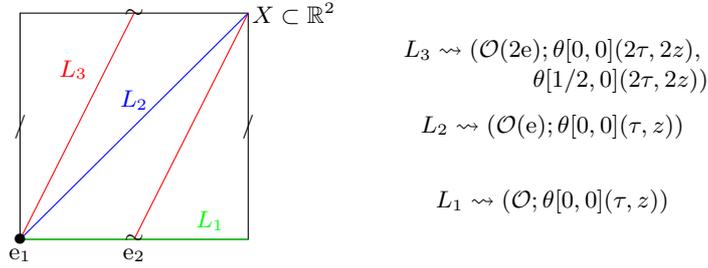
\begin{figure}[ht]
   \centering
  \begin{tikzpicture}
       \draw (0,0) rectangle (3,3);
        \draw[green] (0,0)--(3,0);
        \draw[red] (0,0)--(1.5,3);
        \draw[red] (3,3)--(1.5,0);
        \draw[blue] (0,0)--(3,3);
        \node[red,above,font=\small] at (0.7,2) {$L_3$};
        \node[green,above,font=\small] at (2.5,0) {$L_1$};
        \node[blue,above,font=\small] at (1.5,1.6) {$L_2$};
        \node[font=\small] at (1.5,0) {$\sim$};
        \node[font=\small] at (1.5,3) {$\sim$};
        \node[font=\small] at (0,1.5) {$/$};
        \node[font=\small] at (3,1.5) {$/$};
        \node[font=\small] at (3.6,3) {$X\subset\R^2$};
        \node at (0,0) {$\bullet$};
       \node[below, font=\small] at (0,0) {$\text{e}_1$};
       \node[below, font=\small] at (1.5,0) {$\text{e}_2$};
       \node[font=\small] at (7,0.5) {$L_1 \leadsto (\O; \theta[0,0](\tau,z))$ };
       \node[font=\small] at (7,1.5) {$L_2 \leadsto( \O(\mathrm{e});\theta[0,0](\tau,z))$ };
       \node[font=\small] at (7,2.5) {$L_3 \leadsto (\O(2\mathrm{e}); \theta[0,0](2\tau,2z),$};
       \node[font=\small] at (7.9,2.1) {$\theta[1/2,0](2\tau ,2z))$};
    \end{tikzpicture}
    \caption{In this case, the Lagrangians $L_1, L_2, L_3$ are passing through the origin, hence they have trivial flat connections. Their mirrors as objects of $\mathcal{D}^bCoh(\check{X})$ are respectively the line bundles $\O, \O(\mathrm{e}), \O(2\mathrm{e})$ together with a basis of sections given by theta functions.}
    \label{fig:lagrangians}
\end{figure}

More generally, let $X_q$ be a symplectic torus with $\omega=A dx\wedge dy$, $[b]\in H^2(X_q;\R)/H^2(X_q;\Z)$ and the complexified K\"ahler parameter $q=e^{2\pi i(Ai+b)}$. Let $\check{X}_{\tau}$ be an elliptic curve with complex structure parametrized by $\tau$, then Polishchuk and Zaslow proved the existence of equivalence of categories 

\begin{equation*}
    \mathcal{D}^{b}Coh(\check{X}_q)\cong\mathcal{F}^0(X_q) 
\end{equation*}

with $q=e^{2\pi i\tau}$ and where $\mathcal{F}^0$ denotes a modified Fukaya category (we refer to \cite{PZ00} for details). In the simplified set-up where $\check{X}=\C/\Z\oplus i\Z$ and $X=\R^2/\Z\oplus\Z$ with $\omega=dx\wedge dy$, the mirror functor $\Phi$ is defined as follows: 
\begin{itemize}
    \item it maps line bundles $\check{E}=\O_{\check{X}}(d \mathrm{e})$ of degree $d$ to lines of slope $d$ in $X$ (see Figure \ref{fig:lagrangians});
   \item it maps the theta function $\theta[0,0](\tau,z)$ (that generates global sections of $\O_{\check{X}}(\mathrm{e})$, and $\mathrm{e}$ is the identity for the group structure) to $\{\mathrm{e}_1\}=L_2\cap L_1=L_3\cap L_2$. Then $\theta[0,0](2\tau, 2z)$ is mapped to $\mathrm{e}_1$ too, and $\theta[1/2,0](2\tau,2z)$ is mapped to $\mathrm{e}_2$, where $\{\mathrm{e}_1, \mathrm{e}_2\}=L_1\cap L_3$. 
   \item the compatibility of the composition morphism $m_2$ in $\mathcal{F}$ ($m_2\colon \Hom(L_1,L_2)\otimes \Hom(L_2,L_3)\to\Hom(L_1,L_3) $) with the one in $\mathcal{D}^b$ follows from counting pseudo--holomorphic discs bounding the Lagrangians $L_j$ and the \textit{addition formula} for theta functions.
\end{itemize}

\subsection{Fukaya approach via multi-valued Morse theory:\\
open/closed string conjecture}

In the same paper \cite{fuk05}, Fukaya proposes a mirror symmetry conjecture for a pair $(\check{X},\check{E})$ of a Calabi--Yau manifold $\check{X}$ together with a holomorphic vector bundle $\check{E}$. According to HMS conjecture, the mirror of $(\check{X},\check{E})$ is an object in the Fukaya category, so a pair $(X,L)$ where $L$ is a Lagrangian in $X$ with a unitary local system. Also in this case, Fukaya's main conjecture asserts that mirror symmetry is incarnated at the level of the integral affine manifold $B$ and that multi-valued Morse theory would control both the deformations theory of the pair $(\check{X},\check{E})$ and the holomorphic discs on the mirror $(X,L)$. 

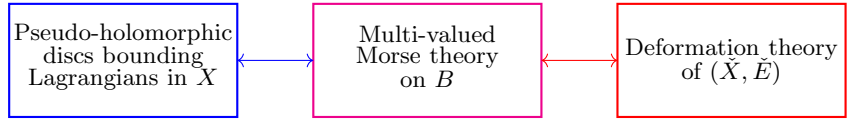
\begin{figure}[ht]
\center
\begin{tikzpicture}
\draw[blue, thick] (0,0) rectangle (3,1.5);
\draw[magenta, thick] (4,0) rectangle (7,1.5);
\draw[red, thick] (8,0) rectangle (11,1.5);
\node[font=\small] at (1.5,1.1) {Pseudo-holomorphic  };
\node[font=\small] at (1.5,0.8) { discs bounding };
\node[font=\small] at (1.5,0.5) {Lagrangians in $X$};
\node[font=\small] at (5.5,1.1) {Multi-valued  };
\node[font=\small] at (5.5,0.8) { Morse theory };
\node[font=\small] at (5.5,0.5) {on $B$};
\node[font=\small] at (9.5,0.9) {Deformation theory  };
\node[font=\small] at (9.5,0.6) { of $(\check{X},\check{E})$ };
\draw[blue, <->]  (3,0.75) -- (4,0.75);
\draw[red, <->]  (7,0.75) -- (8,0.75);
\end{tikzpicture}
\end{figure}

Generalizing the construction of \cite{CLM19}, in \cite{Fan19} we prove that solutions of the Maurer--Cartan equation that governs deformations of holomorphic pairs $(\check{X}=TB_0/\Lambda,\check{E}=\bigoplus_{j=1}^r\O_{\check{X}})$ give rise to consistent scattering diagrams in $B$. Trying to establish a connection between scattering diagrams and pseudo-holomorphic discs counting seems quite difficult in general and working with tropical discs seems rather useful. Indeed, Suen and collaborators introduced tropical Lagrangian multi-sections as combinatorial data which allows reconstructing the mirror bundle (see \cite{suen2019reconstruction,CSM22,Oh-Suen,suen22}). 

A natural question is whether the data of a tropical Lagrangian multi-section can be encoded in scattering diagrams like the one that governs deformations of holomorphic pairs. We leave the discussion to future works. 

\begin{remark}
    In absence of singularities, Abouzaid in \cite{abouzaid14} constructs the mirror of a Lagrangian submanifold $L\subset X$ as a twisted coherent sheaf on $\check{X}$. His construction uses Fukaya's notion of \textit{family of Floer cohomology} and not the multivalued Morse theory approach mentioned above. However, as claimed by Fukaya in \cite{fuk05} (Remark 4.4) the two points of view coincide.
    Along the same line of ideas, recent progress in this direction is due to Yuan \cite{yuan2022family,yuan2021family}. 
\end{remark}

\section{Scattering diagrams}\label{sec:scattering}

Naively, $2$ dimensional scattering diagrams\footnote{There exists a definition of higher dimension scattering diagrams \cite{GS07}.} are collections of lines and rays in $\mathbb{R}^2$ decorated with \textit{automorphisms}. The group structure of scattering diagrams carries the main information, and we are going to discuss it in this section. 

The first aspect to take into consideration is that the automorphisms attached to a ray must belong to a pro-nilpotent Lie group $\hat{G}$. Furthermore, different choices of the group describe different phenomena: for example,
\begin{itemize}
    \item if $\hat{G}$ is the tropical vertex group $\mathbb{V}$ (whose name was introduced by Gross, Pandharipande and Siebert \cite{GPS10}, but the group itself first appeared in \cite{KS-affine}), in \cite{GPS10} the authors show that scattering diagrams in $\mathbb{V}$ compute genus zero Gromov--Witten invariants of $\mathbb{P}^2$ or blow-up of $\mathbb{P}^2$.\footnote{In higher dimensions, scattering diagrams in the so-called \textit{higher dimensional tropical vertex} compute punctured log Gromov-Witten invariants of the log Calabi-Yau variety, as proved in \cite{arguz20}} Moreover, scattering diagrams in the Gross--Siebert program \cite{GS07} (as well as in the GHK construction \cite{GHK}) are a generalization of scattering diagrams in $\mathbb{V}$;
    \item if $\hat{G}$ is the quantum tropical vertex group $\hat{\mathbb{V}}$ (whose name is due to Bousseau \cite{bou_quantum}, but in fact, the group itself was introduced by Kontsevich and Soibelman in \cite{KS-WCF}), scattering diagrams in $\hat{\mathbb{V}}$ compute higher genus log Gromov--Witten invariants for log Calabi--Yau surfaces \cite{bou_quantum}\footnote{In fact, before the work of Bousseau \cite{bou_quantum}, scattering diagrams in $\hat{\mathbb{V}}$ appeared in \cite{FS15} as $q$-deformed diagrams that encode some refined tropical invariants. Bousseau's contribution in \cite{bou_quantum} is to relate these tropical invariants to higher genus curve invariants for log Calabi--Yau surfaces.}. In addition, in \cite{bou_mirror} Bousseau generalizes the GHK construction, promoting scattering diagrams in $\mathbb{V}$ to diagrams in $\hat{\mathbb{V}}$. In fact, he obtained the deformation quantization of the GHK family $\check{\mathcal{X}}\to \mathrm{Spf}\hat{R}$;
    \item if $\hat{G}$ is the extended tropical vertex group $\Tilde{\mathbb{V}}$ introduced in \cite{Fan19}, conjecturally, scattering diagrams in $\Tilde{\mathbb{V}}$ should contribute to the reconstruction problem for holomorphic pairs, in the spirit of \cite{GS07}.    
\end{itemize}

Before giving the definitions of the aforementioned groups, we defined scattering diagrams more generally for a pro-nilpotent Lie group $\hat{G}$. 

Let $\Gamma$ be a lattice, and $\mathfrak{g}=\bigoplus_{\gamma\in\Gamma}\mathfrak{g}_\gamma$ be a $\Gamma$-graded Lie algebra. Assuming the existence of a quadratic form $||-||$ on $\Gamma_{\mathbb{R}}=\Gamma\otimes_{\mathbb{Z}}\mathbb{R}$, we define 
\[\mathfrak{g}_{> n}=\bigoplus_{{\gamma\in\Gamma\setminus 0,\,\, ||\gamma||> n}}\mathfrak{g}_\gamma\]
\[ \mathfrak{g}_{\leq n}=\mathfrak{g}/\mathfrak{g}_{>n}\cong\bigoplus_{\gamma\in\Gamma\setminus 0,\,\, ||\gamma||\leq n}\,\mathfrak{g}_\gamma\]
The completion of $\mathfrak{g}_{\leq n}$ is a pro-nilpotent Lie algebra $\hat{\mathfrak{g}}=\varprojlim_{n}\mathfrak{g}_{\leq n}$ and a pro-nilpotent Lie group is defined via the exponential map
\[ \hat{G}=\exp (\hat{\mathfrak{g}}) \]

and it is a group with the Backer--Campbell--Hausdorff product (BCH)
\[
\exp(g)\cdot\exp (g')=\exp (g+g'+\frac{1}{2}[g,g']...)
\]

Then we define scattering diagrams in $\hat{G}$

\begin{definition}
  Let $\Gamma$ be a rank $2$ lattice and let $\hat{G}$ be a pro-nilpotent Lie group defined from a $\Gamma$-graded Lie algebra. Scattering diagram $\mathfrak{D}$ in $\hat{G}$ is a collection of walls $(\mathfrak{d}_\gamma,\theta_\gamma)$ such that 
  \begin{itemize}
      \item  $\gamma\in\Gamma$
      \item $\mathfrak{d}_\gamma$ is a line (or a ray) $\mathfrak{d}_\gamma=\xi_0+\gamma\mathbb{R}\subset\mathbb{R}^2$ (or $\mathfrak{d}_\gamma=\xi_0+\gamma\mathbb{R}_{\geq 0}\subset\mathbb{R}^2$)
      \item $\theta_\gamma\in \exp( \varprojlim_n \bigoplus_{||k\gamma||\leq n}\mathfrak{g}_{k\gamma})\subset\hat{G} $, in particular $\theta_\gamma$ is the automorphism associated to the line (ray) $\mathfrak{d}_\gamma$
  \end{itemize}
 If $\theta_\gamma=\mathsf{id}$ the wall $(\mathfrak{d}_\gamma,\theta_\gamma)$ does not belong to the diagram $\mathfrak{D}$.  
  
\end{definition}

We remark that $\hat{G}$ being pro-nilpotent is crucial in the previous definition because, on the one hand, the $\Gamma$-grading allows to associate an automorphism to a line of the diagram, and on the other hand, the completion allows to consider products of automorphisms associated to different rays. Notice also that ignoring trivial automorphisms guarantees that products are well defined because modulo $n$ only finitely many rays give meaningful contributions.     

Among all scattering diagrams, the important ones in the reconstruction problem are the so-called \textit{consistent scattering diagrams} $\mathfrak{D}_\infty$ (see Figure \ref{fig:my_label2}). These can be recursively constructed as proved by Kontsevich and Soibelman in \cite{KS-affine}.  

Denote by $\emph{Sing}(\mathfrak{D})$ the singular set of $\mathfrak{D}$: 
\[\emph{Sing}(\mathfrak{D}):=\bigcup_{(\mathfrak{d}_\gamma,\theta_\gamma)\in\mathfrak{D}}\partial\mathfrak{d}_\gamma\cup\bigcup_{\mathfrak{d}_1,\mathfrak{d}_2}\mathfrak{d}_{1}\cap \mathfrak{d}_{2}\] where $\partial \mathfrak{d}=\xi_0$ if $\mathfrak{d}_{\gamma}$ is a ray and zero otherwise. 

There is a notion of order product for the automorphisms associated with each line of a given scattering diagram, and it is defined as follows:

\begin{definition}[Path order product]\label{def:pathorderedprod}
Let $\ell:[0,1]\rightarrow \Gamma_{\mathbb{R}}\setminus\emph{Sing}(\mathfrak{D})$ be a smooth immersion with the starting point that does not lie on a ray of the scattering diagram $\mathfrak{D}$ and such that it intersects transversely the rays of $\mathfrak{D}$. For each power $k>0$, there are times $0< \tau_{1}\leq... \leq\tau_{s} <1$ and rays $\mathfrak{d}_i\in\mathfrak{D}$ such that $\ell(\tau_j)\cap \mathfrak{d}_j\neq 0$. Then, define $\Theta_{\ell,\mathfrak{D}}^k:=\prod_{j=1}^s\theta_{j}$. The path order product is given by: 
\begin{equation}
\Theta_{\gamma,\mathfrak{D}}:=\lim_{k\rightarrow\infty}\Theta_{\ell,\mathfrak{D}}^k
\end{equation}
\end{definition}

\begin{definition}[Consistent scattering diagram]A scattering diagram $\mathfrak{D}_\infty$ is consistent if for any closed path $\ell$ intersecting $\mathfrak{D}_\infty$ generically, $\Theta_{\gamma,\mathfrak{D}_\infty}=\mathsf{id}_{\hat{G}}$. 
\end{definition}
The following theorem by Kontsevich and Soibelman is an existence (ad uniqueness) result of consistent scattering diagrams: 
\begin{theorem}[\cite{KS-affine}]\label{thm:KS}
Let $\mathfrak{D}$ be a scattering diagram with two non-parallel walls. There exists a unique minimal\footnote{The diagram is minimal meaning that we do not consider rays with trivial automorphisms.} scattering diagram $\mathfrak{D}_{\infty}\supseteq\mathfrak{D}$ such that $\mathfrak{D}_{\infty}\setminus\mathfrak{D}$ consists only of rays, and it is {consistent}.
\end{theorem}
The proof of the Theorem is based on a factorization property for elements of $\hat{G}$, namely working order by order in $n$, computing the product of the $\theta_\gamma$ for $||\gamma||\leq n$ and then factorizing the result with respect to different $\gamma$ will give a prescription for new rays. 

We say that the consistent scattering diagram $\mathfrak{D}_\infty$ of Theorem \ref{thm:KS} \textit{saturates} the initial scattering diagram $\mathfrak{D}$.

\subsection{Scattering diagrams in the tropical vertex group}
In this section, we recall the definition of the tropical vertex group following the treatment of \cite{GPS10}. 

The tropical vertex group $\mathbb{V}$ is a subgroup of symplectomorphisms of the formal algebraic torus $\mathbb{T}_{\Lambda}=\Spec\C[\Lambda][\![t]\!]$, $\Lambda\simeq\Z^2$ being a rank two lattice. 
Let $e_1$ and $e_2$ be a basis for $\Lambda\cong\Z^2$. The group ring $\mathbb{C}[\Lambda]$ is the ring of Laurent polynomial in the variable $z^m$, with the convention that $z^{e_1}=:x$ and $z^{e_2}=:y$. We denote by $\Lambda^*=\Hom(\Lambda, \Z)$ be the dual lattice and set $\langle -,-\rangle\colon\Lambda^*\times\Lambda\to\Z$ be the symmetric pairing. Every elements $n\in\Lambda^*$ is associated with a derivation $\partial_n\in\Der(\C[\Lambda])$ whose action is defined as $z^{m}\partial_n(z^{m'}):= z^{m+m'}\langle n,m'\rangle$. We are going to introduce the formal Lie algebra of derivations:

\begin{equation}
\mathfrak{g}:=\big(\mathbb{C}[\Lambda]{\otimes}_{\mathbb{C}}(t)\C[\![ t ]\!]\big)\otimes_{\mathbb{Z}} \Lambda^*
\end{equation}
where the natural Lie bracket on $\mathfrak{g}$ is 
\begin{equation}
[z^m\partial_n,z^{m'}\partial_{n'}]:= z^{m+m'}\partial_{\langle m',n\rangle n'-\langle m, n'\rangle n}.
\end{equation}
In particular $\mathfrak{g}$ has a Lie sub-algebra $\mathfrak{h}\subset\mathfrak{g}$ defined by:
\begin{equation}
\mathfrak{h}:=\big(\bigoplus_{m\in \Lambda\setminus\lbrace 0\rbrace}\C z^m\cdot (t)\otimes_{\Z} m^{\perp}\big)\hat{\otimes}\C[\![ t ]\!]
\end{equation}
where $m^{\perp}\in\Lambda^*$ is identified with the derivation $\partial_n$ for a unique primitive vector $n\in\Lambda^*$ such that $\langle n,m\rangle=0$ and it is positively oriented according with the orientation induced by $\Lambda_{\mathbb{R}}:=\Lambda\otimes_{\Z}\mathbb{R}$. By exponentiating the Lie algebra $\mathfrak{h}$ we get a subgroup of the group of formal automorphisms of the algebraic torus, i.e. the tropical vertex group
\begin{definition}
The tropical vertex group $\mathbb{V}$ is the sub-group of $\Aut_{\C[\![ t ]\!]}\big(\C[\Lambda]\hat{{\otimes}}_{\mathbb{C}}\C[\![ t ]\!]\big)$, such that $\mathbb{V}:=\exp(\mathfrak{h})$. The product on $\mathbb{V}$ is defined by the BCH formula.
\end{definition}
Equivalently, $\mathbb{V}$ can be defined as the group generated by formal one-parameter families of symplectomorphisms of the algebraic torus $\mathbb{T}_{\Lambda}=\Spec\C[\Lambda]$ defined as follows: let $f_{(a,b)}=1+tx^ay^b\cdot g(x^ay^b,t)\in\C[\Lambda]\hat{\otimes}_{\C}\C[\![t]\!]$ for some $(a,b)\in\Lambda$ and for a polynomial $g(x^ay^b,t)\in\C[x^ay^b,t]$, then $\mathbb{V}$ is the group generated by all possible $\theta_{(a,b),f_{(a,b)}}$ defined as 
\begin{equation}
\theta_{(a,b),f_{(a,b)}}(x)=f^{-b}x \qquad \theta_{(a,b),f_{(a,b)}}(y)=f^ay.
\end{equation}

In particular, $\theta_{(a,b),f_{(a,b)}}^*\omega=\omega$ with respect to the holomorphic symplectic form on $\omega=\frac{dx}{x}\wedge \frac{dy}{y}$.

The two definitions are equivalent as one can easily check that  
\[\theta_{(a,b),f_{(a,b)}}=\exp\left(\log (f_{(a,b)}) \partial_{(-b,a)}\right)\in\mathbb{V}.\] 

We can now state the definition of scattering diagrams: 

\begin{definition}[Scattering diagram]\label{def:scattering}
A scattering diagram $\mathfrak{D}$ (in $\mathbb{V}$) is a collection of \textit{walls} $\textsf{w}_i=(\mathfrak{d}_i,\theta_i)$, where 
\begin{itemize}
   \item $\mathfrak{d}_i$ can be either a \textit{line} through $\xi_0$, i.e. $\mathfrak{d}_i=\xi_0-m_i\mathbb{R}$ or a \textit{ray} (half line) $\mathfrak{d}_i=\xi_0-m_i\mathbb{R}_{\geq 0}$,
    \item $\theta_i\in\mathbb{V}$ is such that $\log(\theta_i)=\sum_{j,k}a_{jk}t^j z^{km_i}\partial_{n_i}\in\mathfrak{h}$. 
\end{itemize}
Moreover for any $N>0$ there are finitely many $\theta_i$ such that $\theta_i\not\equiv 1$ mod $t^N$.
\end{definition}

\begin{remark}
Notice that $\mathfrak{h}$ is both $\Lambda$-granded and $t$-graded. However, in the definition of $\mathbb{V}$ we choose the completion in the $t$ parameter and not with respect to $\Lambda$ as in the general definition. That explains why we add the latter assumption in the definition of scattering diagrams in $\mathbb{V}$, as otherwise, we would not be able to define products in $\mathbb{V}$.
\end{remark}


In \cite{CLM19}, the authors prove that the Lie algebra $\mathfrak{h}$ embeds to the degree zero part of the Kodaira--Spencer DGLA which governs infinitesimal deformations of $\check{X}:=TB/\Lambda\to B$ with $B$ being a tropical affine smooth $2$-dimensional manifold (as the toy model example for SYZ).\footnote{In \cite{CLM19}, they do not restrict to $2$ dimensional manifolds and they work with $n$-dimensional scattering diagrams.} In particular, the group $\mathbb{V}$ acts as the Gauge group on the solutions of the Maurer--Cartan equation. Hence, they show that consistent scattering diagrams in $\mathbb{V}$ can be constructed by studying the asymptotics of formal solutions of the Maurer--Cartan equation. In \cite{chan2020tropical}, they also explain how consistent scattering diagrams in $\mathbb{V}$ encode tropical discs counting. The latter should correspond to count pseudo--holomorphic discs on $X$, so summarizing, we have the following picture
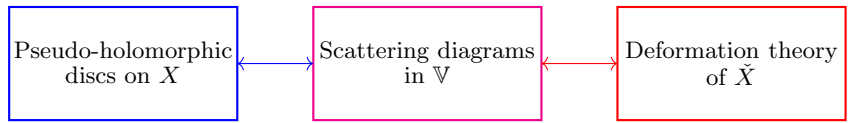
\begin{figure}[ht]
\center
\begin{tikzpicture}
\draw[blue, thick] (0,0) rectangle (3,1.5);
\draw[magenta, thick] (4,0) rectangle (7,1.5);
\draw[red, thick] (8,0) rectangle (11,1.5);
\node[font=\small] at (1.5,0.9) {Pseudo-holomorphic  };
\node[font=\small] at (1.5,0.6) { discs on $X$ };
\node[font=\small] at (5.5,0.9) {Scattering diagrams };
\node[font=\small] at (5.5,0.6) { in $\mathbb{V}$};
\node[font=\small] at (9.5,0.9) {Deformation theory  };
\node[font=\small] at (9.5,0.6) { of $\check{X}$ };
\draw[blue, <->]  (3,0.75) -- (4,0.75);
\draw[red, <->]  (7,0.75) -- (8,0.75);
\end{tikzpicture}
\end{figure}

\begin{remark}
The group $\mathbb{V}$ do also appear in \cite{KS-WCF} (see Example 4 in Section 2.3, when $\Gamma$ is $2$-dimensional) as the pronilpotent Lie group of the \textit{torus Lie algebra} (which up to a sign, is equivalent to $\mathfrak{h}$). In particular, one can consider automorphisms in $\mathbb{V}$ which encode DT-invariants and define the so-called wall-crossing structure. The latter appears in \cite{GMN}, where they allow to compute the BPS spectrum of certain theories. Finally, consistent scattering diagrams in $\mathbb{V}$ are studied in \cite{Bri} in relation to quiver representations.     
\end{remark}

\subsection{Scattering diagrams in the quantum tropical vertex}

In \cite{KS-WCF}, the authors first considered the quantization of the Poisson algebra of functions of a torus $\mathbb{T}_{\Lambda}=\Spec \mathbb{C}[\Lambda]$, $\Lambda\cong \mathbb{Z}^2$ in the context of wall--crossing formulas for DT invariants on 3d Calabi--Yau categories. We briefly recall the construction following \cite{bou_quantum}. Let $\mathcal{O}_{\mathbb{T}}$ be the ring of functions of $\mathbb{T}_{\Lambda}$ with $\Gamma(\mathcal{O}_{\mathbb{T}})=\bigoplus_{m\in\Lambda}\C\, z^m$ and 

\begin{equation*}
    z^m \cdot z^{m'}= z^{m+m'}. 
\end{equation*}

Then the Poisson bracket is 
\begin{equation*}
    \lbrace z^m, z^{m'}\rbrace = \langle m,m'\rangle_D z^{m+m'}
\end{equation*}
with $\langle -,-\rangle_D$ being an unimodular antisymmetric bilinear form on $\Lambda$. The ring of functions of the \textit{quantum torus} $\hat{\mathbb{T}}$ is defined as $\Gamma(\mathcal{O}_{\hat{\mathbb{T}}}):=\bigoplus_{m\in\Lambda}\C[q^{\pm 1/2}]\, \hat{z}^m$ and 

\begin{equation*}
    \hat{z}^m\cdot \hat{z}^{m'}= q^{\frac{1}{2}\langle m, m'\rangle_D} \hat{z}^{m+m'}. 
\end{equation*}


Analogously to the tropical vertex group, $\hat{\mathbb{V}}$ is the group of formal automorphisms of $\hat{\mathbb{T}}$: let us first consider the pro-nilpotent, $\Lambda$-graded Lie algebra

\begin{equation}
    \hat{\mathfrak{h}}:=(t)\left(\bigoplus_{m\in\Lambda\setminus\{0\}} \C[q^{\pm\frac{1}{2}}] \hat{z}^m\right)  \hat{\otimes}_{\C}\C[\![t]\!]
\end{equation}

whose Lie bracket simply reads $[\hat{z}^m, \hat{z}^{m'}]=(q^{\frac{1}{2}\langle m, m'\rangle_D}-q^{-\frac{1}{2}\langle m,m'\rangle_D}) \hat{z}^{m+m'} $. 

\begin{definition}
    The quantum tropical vertex group $\hat{\mathbb{V}}$ is the subgroup of the formal automorphisms $\Aut_{\C[\![t]\!]}\left(\C[q^{\pm\frac{1}{2}}][\Lambda]\hat{\otimes}_{\C}\C[\![t]\!]\right)$ such that $\hat{\mathbb{V}}=\exp(\hat{\mathfrak{h}})$. The product is defined by the BCH formula. 
\end{definition}
In particular, elements of $\hat{\mathbb{V}}$ can be defined as 

\begin{align*}
    \hat{\theta}_{m}&:=\mathrm{Ad}_{\hat{H}_m}=\exp(\hat{H}_m) (-)\exp(-\hat{H}_m)
    \end{align*}
for $\hat{H}_m :=\sum_{k,j=1}^{\infty} a_{k,j}(q)t^j \hat{z}^{mk}\in\hat{\mathfrak{h}}$.


\subsection{Scattering diagrams in the extended tropical vertex group}

In \cite{Fan19}, the author studied the relationship between deformations of holomorphic pairs and scattering diagrams. The former consists of a rank $r$ holomorphic vector bundle $\check{E}$ over a semi-flat manifold $\check{X}:=TB/\Lambda\to B$, namely the toy model torus fibration over a $2$ dimensional smooth affine tropical manifold $B$. The construction generalizes the one of Chan, Conan Leung and Ma \cite{CLM19} and we introduce the extended tropical vertex group $\tilde{\mathbb{V}}$ where the scattering diagrams are defined.

The group $\tilde{\mathbb{V}}$ is defined by taking the exponential of the Lie algebra $\tilde{\mathfrak{h}}$ which is a $\mathfrak{gl}$-extension of $\mathfrak{h}$
\begin{equation}
\tilde{\mathfrak{h}}:=\left(\bigoplus_{m\in \Lambda\setminus\lbrace 0\rbrace}\C z^{m}(t)\otimes_{\Z}\left(\mathfrak{gl}(r,\mathbb{C})\oplus  m^{\perp}\right)\right)\hat{\otimes}_{\C}\C[\![t]\!] 
\end{equation}
with Lie bracket $[-,- ]_{\sim}\colon \tilde{\mathfrak{h}}\times \tilde{\mathfrak{h}}\to \tilde{\mathfrak{h}}$ is 
\begin{multline}\label{eq:Liebracket}
[(A,\partial_n)z^m , (A',\partial_{n'})z^{m'} ]_{\sim}:=([A,A']_{\mathfrak{gl}}z^{m+m'}+A'\langle m',n\rangle z^{m+m'}- A\langle m,n'\rangle z^{m+m'},\\
[z^m\partial_n,z^{m'}\partial_{n'}]_{\mathfrak{h}} ).
\end{multline}
The definition of the Lie algebra $\tilde{\mathfrak{h}}$ comes from the Kodaira Spencer DGLA which governs infinitesimal deformations of holomorphic pairs \cite{def_pair}. We refer to \cite{Fan19} for more details about the construction. 

\begin{definition}
    The extended tropical vertex group $\tilde{\mathbb{V}}$ is defined as the pro-nilpotent group $\tilde{\mathbb{V}}=\exp(\tilde{\mathfrak{h}})$. The product on $\tilde{\mathbb{V}}$ is defined by the BCH formula.
\end{definition}

Following Fukaya's approach to mirror symmetry, we expect consistent scattering diagrams in $\tilde{\mathbb{V}}$ to govern the reconstruction of the mirror pair $(L_1,...,L_r,X)$ conjecturally given by the mirror manifold $X$ and a collection of $r$ Lagrangians $L_j$ in $X$ (as previously discussed in Section \ref{HMS}). In addition, we expect that holomorphic discs bounding $L_1,...,L_r$ and the fibers of $X$ are encoded in the commutators of some elements in $\Tilde{\mathbb{V}}$. In \cite{Fan20} we consider some \textit{ad-hoc} scattering diagrams in $\tilde{\mathbb{V}}$ and computing commutators we recover genus zero Gromov--Witten invariants of blow-up of $\mathbb{P}^2$ maximally tangent to its boundary divisor as well as genus zero Gromov--Witten invariants relative to the toric boundary divisor. Alternatively, we consider $2d$-$4d$ scattering diagrams in $\Tilde{\mathbb{V}}$ which seems more naturally related to the reconstruction problem for the mirror pair.\footnote{Recall that pure $4d$ scattering diagrams are defined in the tropical vertex group $\mathbb{V}$, hence when the $2d$ component is trivial, we end up studying the reconstruction problem of $X$ itself.} 
For instance, let $\mathfrak{D}$ be the scattering diagram consisting of two initial non-parallel walls $\mathsf{w}_i=(\mathfrak{d}_i,\theta_i)$ $i=1,2$

\begin{center}
\begin{tabular}{c c c c}
    $\mathsf{w}_1$:  & $\mathfrak{d}_1=m_1\R$ &  $\log\theta_1=(-tE_{ij}z^{m_1},0)$ & \textit{ type } 2d\\  
    $\mathsf{w}_2$: & $\mathfrak{d}_2=m_2\R$ &
   $\log\theta_2=(0, \sum_{k\geq 0}\frac{z^{km_2}}{k}t^k\partial_{n_2})$& \textit{ type } 4d
\end{tabular}
\end{center}

where $E_{ij}\in\mathfrak{gl}(r,\C)$ is an elementary matrix with non zero entries in position $i<j$, and $m_1,m_2\in\Lambda$. Notice that \textit{type} $2d$ consists only of a \textit{matrix component} while \textit{type} $4d$ consists of a \textit{derivative component} (namely an element of the Lie algebra $\mathfrak{h}$). Indeed the labels $2d$ and $4d$ refer to the formalism of $2d$-$4d$ wall--crossing formula introduced by Gaiotto--Moore--Neitzke \cite{GMN2d4d}. In \cite{Fan19} (Example 1) we compute the consistent scattering diagram $\mathfrak{D}_\infty$ that saturates $\mathfrak{D}$, namely $\mathfrak{D}_{\infty}=\mathfrak{D}\cup\{\mathsf{w}_3\}$  
\begin{center}
    \begin{tabular}{c c c c}
    $\mathsf{w}_3$:  & $\mathfrak{d}_3=(m_1+m_2)\R_{\geq 0}$ &  $\log\theta_3=(t^2E_{ij} \langle m_1,n_2\rangle z^{m_1+m_2},0)$ & \textit{ type } 2d.
    \end{tabular}
\end{center}

We may interpret this result by saying that in presence of a singular torus fiber $X_s$ and of two ramified \textit{horizontal} Lagrangians $L_i,L_j$, the scattering diagram $\mathfrak{D}_\infty$ encodes the number of holomorphic discs bounding $L_i,L_j$ and $X_s$. 

 Furthermore, in \cite{Fan19}, we prove that consistent scattering diagrams in $\tilde{\mathbb{V}}$ are equivalent to $2d$-$4d$ wall--crossing formulas (WCF) which govern BPS state counting in $2d$-$4d$ coupled systems. 


Summing up, the results of \cite{Fan19,Fan20} propose a new approach to the original Fukaya's open--closed string conjecture, with scattering diagrams playing the role of gradient flow trajectories; but we do not know yet an enumerative geometric interpretation.  

\begin{figure}[ht]
\center
\begin{tikzpicture}
\draw[blue, thick] (0,0) rectangle (3,1.5);
\draw[magenta, thick] (4,0) rectangle (7,1.5);
\draw[red, thick] (8,0) rectangle (11,1.5);
\node[font=\small] at (1.5,1.1) {Pseudo-holomorphic  };
\node[font=\small] at (1.5,0.8) { discs bounding };
\node[font=\small] at (1.5,0.5) {Lagrangians in $X$};
\node[font=\small] at (5.5,0.9) {Scattering diagrams };
\node[font=\small] at (5.5,0.6) { in $\tilde{\mathbb{V}}$};
\node[font=\small] at (9.5,0.9) {Deformation theory  };
\node[font=\small] at (9.5,0.6) { of $(\check{X},\check{
E})$ };
\draw[blue, <->, dashed]  (3,0.75) -- (4,0.75);
\draw[red, <->]  (7,0.75) -- (8,0.75);
\end{tikzpicture}
\end{figure}

\printbibliography[heading=bibintoc]

@book{Dbranes,
  title={Dirichlet branes and mirror symmetry},
  author={Aspinwall, Paul and Bridgeland, Tom and Craw, Alastair},
  volume={4},
  year={2009},
  publisher={American Mathematical Soc.}
}

@article{PZ00,
  title={Categorical mirror symmetry: the elliptic curve},
  author={Polishchuk, Alexander and Zaslow, Eric},
  journal={arXiv math/9801119},
  year={1998}
}

@article{KS-HMS,
  title={Homological mirror symmetry and torus fibrations},
  author={Kontsevich, Maxim and Soibelman, Yan},
  journal={arXiv math/0011041},
  year={2000}
}

@incollection{KS-affine,
  title={{A}ffine structures and non-{A}rchimedean analytic spaces},
  author={Kontsevich, Maxim and Soibelman, Yan},
  booktitle={The unity of mathematics},
  pages={321--385},
  year={2006},
  publisher={Springer}
}

@article{SYZ,
  title={{M}irror symmetry is {T}-duality},
  author={Strominger, Andrew and Yau, Shing-Tung and Zaslow, Eric},
  journal={Nuclear Physics B},
  volume={479},
  number={1-2},
  pages={243--259},
  year={1996},
  publisher={Elsevier}
}

@book{CoxKatz,
  title={Mirror symmetry and algebraic geometry},
  author={Cox, David A and Katz, Sheldon},
  volume={68},
  year={1999},
  publisher={American Mathematical Society Providence, RI}
}

@article{GS1,
  title={Mirror symmetry via logarithmic degeneration data {I}},
  author={Gross, Mark and Siebert, Bernd},
  journal={Journal of Differential Geometry},
  volume={72},
  number={2},
  pages={169--338},
  year={2006},
  publisher={Lehigh University}
}

@article{GS07,
  title={From real affine geometry to complex geometry},
  author={Gross, Mark and Siebert, Bernd},
  journal={Annals of mathematics},
  pages={1301--1428},
  year={2011},
  publisher={JSTOR}
}

@article{GSinvitation,
  title={An invitation to toric degenerations},
  author={Gross, Mark and Siebert, Bernd},
  journal={arXiv:0808.2749},
  year={2008}
}

@article{GS2,
  title={Mirror symmetry via logarithmic degeneration data, {II}},
  author={Gross, Mark and Siebert, Bernd},
  journal={Journal of Algebraic Geometry},
  volume={19},
  number={4},
  pages={679--780},
  year={2010}
}

@article{GSlocal,
  title={Local mirror symmetry in the tropics},
  author={Gross, Mark and Siebert, Bernd},
  journal={arXiv:1404.3585},
  year={2014}
}

@incollection{GS-review,
  title={Enumerative aspects of the {G}ross-{S}iebert program},
  author={van Garrel, Michel and Overholser, D Peter and Ruddat, Helge},
  booktitle={Calabi-Yau varieties: arithmetic, geometry and physics},
  pages={337--420},
  year={2015},
  publisher={Springer}
}

@article{GSintrinsic,
  title={Intrinsic mirror symmetry and punctured {G}romov-{W}itten invariants},
  author={Gross, Mark and Siebert, Bernd},
  journal={Algebraic geometry: Salt Lake City 2015},
  volume={97},
  pages={199--230},
  year={2018}
}

@article{GS19,
  title={Intrinsic mirror symmetry},
  author={Gross, Mark and Siebert, Bernd},
  journal={arXiv:1909.07649},
  year={2019}
}

@article{GS22,
  title={The canonical wall structure and intrinsic mirror symmetry},
  author={Gross, Mark and Siebert, Bernd},
  journal={Inventiones mathematicae},
  volume={229},
  number={3},
  pages={1101--1202},
  year={2022},
  publisher={Springer}
}

@article{CLM19,
  title={Scattering diagrams from asymptotic analysis on {M}aurer--{C}artan equations},
  author={Chan, Kwokwai and Leung, Naichung Conan and Ma, Ziming Nikolas},
  journal={Journal of the European Mathematical Society},
  volume={24},
  number={3},
  pages={773--849},
  year={2021}
}

@article{Fan19,
  title={Deformations of holomorphic pairs and 2d-4d wall-crossing},
  author={Fantini, Veronica},
  journal={Advances in Theoretical and Mathematical Physics},
  volume={6},
number={26},
pages={1705--1769},
  year={2022},
publisher={International Press}
}

@inproceedings{fuk05,
  title={Asymptotic analysis and mirror symmetry},
  author={Fukaya, Kenji},
  booktitle={Graphs and Patterns in Mathematics and Theoretical Physics: Proceedings of the Conference on Graphs and Patterns in Mathematics and Theoretical Physics, Dedicated to Dennis Sullivan's 60th Birthday, June 14-21, 2001, Stony Brook University, Stony Brook, NY},
  volume={73},
  pages={205},
  year={2005},
  organization={American Mathematical Soc.}
}

@article{bou_quantum,
  title={The quantum tropical vertex},
  author={Bousseau, Pierrick},
  journal={Geometry \& Topology},
  volume={24},
  number={3},
  pages={1297--1379},
  year={2020},
  publisher={Mathematical Sciences Publishers}
}

@article{bou_mirror,
  title={Quantum mirrors of log {C}alabi--{Y}au surfaces and higher-genus curve counting},
  author={Bousseau, Pierrick},
  journal={Compositio Mathematica},
  volume={156},
  number={2},
  pages={360--411},
  year={2020},
  publisher={London Mathematical Society}
}

@article{arguz20,
  title={The higher dimensional tropical vertex},
  author={Arg{\"u}z, H{\"u}lya and Gross, Mark},
  journal={arXiv:2007.08347},
  year={2020}
}

@article{GPS10,
  title={The tropical vertex},
  author={Gross, Mark and Pandharipande, Rahul and Siebert, Bernd},
  journal={Duke Mathematical Journal},
  volume={153},
  number={2},
  pages={297--362},
  year={2010},
  publisher={Duke University Press}
}

@article{FS15,
  title={Block--{G}{\"o}ttsche invariants from wall-crossing},
  author={Filippini, Sara Angela and Stoppa, Jacopo},
  journal={Compositio Mathematica},
  volume={151},
  number={8},
  pages={1543--1567},
  year={2015},
  publisher={London Mathematical Society}
}

@article{CSM22,
  title={Tropical {L}agrangian multi-sections and smoothing of locally free sheaves over degenerate {C}alabi--{Y}au surfaces},
  author={Chan, Kwokwai and Ma, Ziming Nikolas and Suen, Yat-Hin},
  journal={Advances in Mathematics},
  volume={401},
  pages={108280},
  year={2022},
  publisher={Elsevier}
}

@article{suen22,
  title={Tropical Lagrangian multi-sections and tropical locally free sheaves},
  author={Suen, Yat-Hin},
  journal={arXiv:2203.02162},
  year={2022}
}

@article{Oh-Suen,
  title={Lagrangian multi-sections and their toric equivariant mirror},
  author={Oh, Yong-Geun and Suen, Yat-Hin},
  journal={arXiv:2211.12191},
  year={2022}
}

@article{suen2019reconstruction,
  title={Reconstruction of $T_{\mathbb{P}^2}$ via tropical {L}agrangian multi-section},
  author={Suen, Yat-Hin},
  journal={arXiv:1904.12449},
  year={2019}
}

@article{KS-WCF,
  title={Wall-crossing structures in {D}onaldson--{T}homas invariants, integrable systems and mirror symmetry},
  author={Kontsevich, Maxim and Soibelman, Yan},
  booktitle={Homological mirror symmetry and tropical geometry},
  pages={197--308},
  year={2014},
  publisher={Springer}
}

@article{auroux2009special,
  title={Special Lagrangian fibrations, wall-crossing, and mirror symmetry},
  author={Auroux, Denis},
  journal={arXiv:0902.1595},
  year={2009}
}

@article{Fan20,
  title={The extended tropical vertex group},
  author={Fantini, Veronica},
  journal={arXiv:2012.05069},
  year={2020}
}

@article{COGP,
    author = "Candelas, Philip and De La Ossa, Xenia C. and Green, Paul S. and Parkes, Linda",
    editor = "Yau, Shing-Tung",
    title = "{A Pair of {C}alabi-{Y}au manifolds as an exactly soluble superconformal theory}",
    reportNumber = "UTTG-25-90",
    doi = "10.1016/0550-3213(91)90292-6",
    journal = "Nucl. Phys. B",
    volume = "359",
    pages = "21--74",
    year = "1991"
}

@article{def_pair,
  title={A differential-geometric approach to deformations of pairs $(X, E)$},
  author={Chan, Kwokwai and Suen, Yat-Hin},
  journal={Complex Manifolds},
  volume={3},
  number={1},
  year={2016},
  publisher={De Gruyter Open Access}
}

@article{GHK,
  title={Mirror symmetry for log {C}alabi-{Y}au surfaces {I}},
  author={Gross, Mark and Hacking, Paul and Keel, Sean},
  journal={Publications Math{\'e}matiques de l'IHES},
  volume={122},
  number={1},
  pages={65--168},
  year={2015},
  publisher={Springer}
}

@article{Gross-P2,
  title={Mirror symmetry for $P^2$ and tropical geometry},
  author={Gross, Mark},
  journal={Advances in Mathematics},
  volume={224},
  number={1},
  pages={169--245},
  year={2010},
  publisher={Elsevier}
}

@misc{morrison1992essays,
  title={Essays on {M}irror manifolds},
  author={Morrison, D},
  year={1992},
  publisher={Intenational Press}
}

@article{barrott,
  title={Explicit equations for mirror families to log {C}alabi-{Y}au surfaces},
  author={Barrott, Lawrence Jack},
  journal={arXiv:1810.08356},
  year={2018}
}

@article{chan2022smoothing,
  title={Smoothing, scattering, and a conjecture of {F}ukaya},
  author={Chan, Kwokwai and Leung, Naichung Conan and Ma, Ziming Nikolas},
  journal={arXiv:2205.09926},
  year={2022}
}

@article{chan2020tropical,
  title={Tropical counting from asymptotic analysis on {M}aurer-{C}artan equations},
  author={Chan, Kwokwai and Ma, Ziming Nikolas},
  journal={Transactions of the American Mathematical Society},
  volume={373},
  number={9},
  pages={6411--6450},
  year={2020}
}

@book{seidel2008fukaya,
  title={Fukaya categories and {P}icard-{L}efschetz theory},
  author={Seidel, Paul},
  volume={10},
  year={2008},
  publisher={European Mathematical Society}
}

@inproceedings{K95-HMS,
  title={Homological algebra of mirror symmetry},
  author={Kontsevich, Maxim},
  booktitle={Proceedings of the International Congress of Mathematicians: August 3--11, 1994 Z{\"u}rich, Switzerland},
  pages={120--139},
  year={1995},
  organization={Springer}
}

@book{gross2012calabi,
  title={Calabi-{Y}au manifolds and related geometries: lectures at a summer school in {N}ordfjordeid, {N}orway, {J}une 2001},
  author={Gross, Mark and Huybrechts, Daniel and Joyce, Dominic},
  year={2012},
  publisher={Springer Science \& Business Media}
}

@book{MS,
  title={Mirror symmetry},
  author={Hori, Kentaro},
  volume={1},
  year={2003},
  publisher={American Mathematical Soc.}
}

@article{GMN,
  title={Four-dimensional wall-crossing via three-dimensional field theory},
  author={Gaiotto, Davide and Moore, Gregory W and Neitzke, Andrew},
  journal={Communications in Mathematical Physics},
  volume={299},
  number={1},
  pages={163--224},
  year={2010},
  publisher={Springer}
}

@article{Bri,
  title={Scattering diagrams, {H}all algebras and stability conditions},
  author={Bridgeland, Tom},
  journal={arXiv:1603.00416},
  year={2016}
}

@article{GMN2d4d,
  title={Wall-crossing in coupled 2d-4d systems},
  author={Gaiotto, Davide and Moore, Gregory W and Neitzke, Andrew},
  journal={Journal of High Energy Physics},
  volume={2012},
  number={12},
  pages={1--169},
  year={2012},
  publisher={Springer}
}

@article{candelas1,
  title={Calabi-{Y}au manifolds in weighted $P^4$},
  author={Candelas, Philip and Lynker, Monika and Schimmrigk, Rolf},
  journal={Nuclear Physics B},
  volume={341},
  number={2},
  pages={383--402},
  year={1990},
  publisher={Elsevier}
}

@article{greene1990duality,
  title={Duality in {C}alabi-{Y}au moduli space},
  author={Greene, Brian R and Plesser, M Ronen},
  journal={Nuclear Physics B},
  volume={338},
  number={1},
  pages={15--37},
  year={1990},
  publisher={Elsevier}
}

@article{arguz-higher-logCY,
  title={Equations of mirrors to log {C}alabi--{Y}au pairs via the heart of canonical wall structures},
  author={Arg{\"u}z, H{\"u}lya},
  journal={arXiv:2109.08664},
  year={2021}
}

@article{abramovich2020punctured,
  title={Punctured logarithmic maps},
  author={Abramovich, Dan and Chen, Qile and Gross, Mark and Siebert, Bernd},
  journal={arXiv:2009.07720},
  year={2020}
}

@article{tu14,
  title={On the reconstruction problem in mirror symmetry},
  author={Tu, Junwu},
  journal={Advances in Mathematics},
  volume={256},
  pages={449--478},
  year={2014},
  publisher={Elsevier}
}

@article{abouzaid14,
  title={Family {F}loer cohomology and mirror symmetry},
  author={Abouzaid, Mohammed},
  journal={arXiv:1404.2659},
  year={2014}
}

@inproceedings{FFR,
  title={Smoothing toroidal crossing spaces},
  author={Felten, Simon and Filip, Matej and Ruddat, Helge},
  booktitle={Forum of Mathematics, Pi},
  volume={9},
  pages={e7},
  year={2021},
  organization={Cambridge University Press}
}

@article{chan21--smoothing,
  title={Geometry of the {M}aurer-{C}artan equation near degenerate {C}alabi-{Y}au varieties},
  author={Chan, Kwokwai and Leung, Naichung Conan and Ma, Ziming Nikolas},
  journal={arXiv:1902.11174},
  year={2019}
}

@article{CJL21,
  title={Special Lagrangian submanifolds of log {C}alabi--{Y}au manifolds},
  author={Collins, Tristan C and Jacob, Adam and Lin, Yu-Shen},
  journal={Duke Mathematical Journal},
  volume={170},
  number={7},
  pages={1291--1375},
  year={2021},
  publisher={Duke University Press}
}

@article{CJL20,
  title={The {SYZ} mirror symmetry conjecture for del {P}ezzo surfaces and rational elliptic surfaces},
  author={Collins, Tristan C and Jacob, Adam and Lin, Yu-Shen},
  journal={arXiv:2012.05416},
  year={2020}
}

@article{yuan2022family,
  title={Family {F}loer mirror space for local {SYZ} singularities},
  author={Yuan, Hang},
  journal={arXiv:2206.04652},
  year={2022}
}

@phdthesis{yuan2021family,
  title={Family {F}loer program and non-archimedean {SYZ} mirror construction},
  author={Yuan, Hang},
  year={2021},
  school={State University of New York at Stony Brook}
}

@article{joyce2000singularities,
  title={Singularities of special {L}agrangian fibrations and the {SYZ} {C}onjecture},
  author={Joyce, Dominic},
  journal={arXiv math/0011179},
  year={2000}
}

@article{groman2022closed,
  title={Closed string mirrors of symplectic cluster manifolds},
  author={Groman, Yoel and Varolgunes, Umut},
  journal={arXiv:2211.07523},
  year={2022}
}

@article{groman2021locality,
  title={Locality of relative symplectic cohomology for complete embeddings},
  author={Groman, Yoel and Varolgunes, Umut},
  journal={arXiv:2110.08891},
  year={2021}
}

@article{KY19,
  title={The {F}robenius structure theorem for affine log {C}alabi-{Y}au varieties containing a torus},
  author={Keel, Sean and Yu, Tony Yue},
  journal={arXiv:1908.09861},
  year={2019}
}

@article{gross-wilson,
  title={Large complex structure limits of {K}3 surfaces},
  author={Gross, Mark and Wilson, Pelham MH},
  journal={Journal of Differential Geometry},
  volume={55},
  number={3},
  pages={475--546},
  year={2000},
  publisher={Lehigh University}
}

\end{document}